\documentclass[a4paper]{article}
\title{A Koopman Operator Tutorial with Orthogonal Polynomials}
\author{Simone Servadio\thanks{Postdoctoral Associate, Department of Aeronautics and Astronautics, Massachusetts Institute of Technology, simoserv@mit.edu}, David Arnas\thanks{Assistant Professor, School or Aeronautics and Astronautics, Purdue University, darnas@purdue.edu}, and Richard Linares\thanks{Boeing Assistant Professor, Department of Aeronautics and Astronautics, Massachusetts Institute of Technology, linaresr@mit.edu}}
\date{\today}

\usepackage[english]{babel}
\usepackage[utf8]{inputenc}
\usepackage{amsmath}
\usepackage{graphicx}
\usepackage[colorinlistoftodos]{todonotes}
\usepackage{float}
\usepackage[numbered]{mcode}
\usepackage[english]{babel}
\usepackage{blindtext}
\usepackage{xfrac}    
\usepackage{algorithm,algorithmic}
\usepackage{subfigure}
\usepackage{psfrag}
\usepackage{bm}
\usepackage{graphicx}
\usepackage{amsmath}
\usepackage{amsfonts}
\usepackage[version=4]{mhchem}
\usepackage{siunitx}
\usepackage{longtable,tabularx}
\newcommand{\mbf}{\mathbf}

\begin{document}
\maketitle

\begin{abstract}
The Koopman Operator (KO) offers a promising alternative methodology to solve ordinary differential equations analytically. The solution of the dynamical system is analyzed in terms of observables, which are expressed as a linear combination of the eigenfunctions of the system. Coefficients are evaluated via the Galerkin method, using Legendre polynomials as a set of orthogonal basis functions. This tutorial provides a detailed analysis of the Koopman theory, followed by a rigorous explanation of the KO implementation in a computer environment, where a line-by-line description of a MATLAB code solves the Duffing oscillator application. 
\end{abstract}

\tableofcontents

\section{Introduction}
The solution to ordinary differential equations (ODEs) can be obtained  both numerically and analytically. The Koopman Operator (KO) approach analyzes the spectral behaviour of the system and provides a deep understanding of its dynamics. The time evolution of the state of the system is provided as a linear combination of a set of well-defined functions instead of performing a step-by-step propagation as done in numerical methods such as Runge-Kutta techniques.

The goal of the proposed KO-based method is reformulating nonlinear problems into a linear framework that can be solved with available linear techniques. Motivated by Koopman~\cite{koopman1931hamiltonian} and Von Neumann~\cite{neumann1932operatorenmethode}, a new  ``Heisenberg picture” was introduced in classical dynamics. In the work by Koopman, it was observed that a KO for a Hamiltonian system is unitary in an $L_2$ Hilbert space. Later Von Neumann was able to make a connection between the spectrum of the KO and that of ergodicity of classical dynamical systems. The key result was demonstrating that there exists an infinite-dimensional linear operator, given by $\mathcal{K}$, that evolves all observation functions $g({\bf x})$ of the state, ${\bf x}$, for any nonlinear system. The evolution of these observables and the KO is defined by the chain rule, which is a linear operator, giving KO its linear properties. 

The linearity of the KO is very appealing, but this benefit is contrasted with the fact that it is infinite-dimensional. However, this issue can be overcome by capturing the evolution on a finite subspace spanned by a finite set of basis functions instead of capturing the evolution of all measurement functions in a Hilbert space. In effect, this is a truncation of the KO to a finite subspace. Additionally, a Koopman invariant subspace is spanned by a set of eigenfunctions of the Koopman operator. A Koopman eigenfunction, $\phi_i({\bf x})$, corresponding to eigenvalue $\lambda_i$ is invariant under the Koopman operator (but for a normalized constant). As such, the evolution of the Koopman eigenfunctions can be expressed as $\frac{d}{dt}\phi_i({\bf x})=\lambda_i\phi_i({\bf x})$. 

This tutorial evaluates the Koopman matrix analytically via the Galerkin methodology~\cite{librationkoopman}. A set of orthogonal polynomials is selected to describe the state-space of the dynamics and to represent the eigenfunctions of the system. The Legendre polynomials have been selected for this tutorial due to their properties, which ease the evaluation of integrals onto a selected support~\cite{zonalkoopman}. Other approaches, based on numerical data, have been previously implemented to evaluate the KO matrix, such as the Extended Dynamic Mode Decomposition (EMDM). It has been proven that the computation of the eigenvalues and modes of the system through the Galerkin method is more precise than the data-based approach \cite{servadio2021koopman}.

\section{Koopman Operator Theory}

A classical definition of nonlinear dynamical systems is given by the initial value problem, which can be represented by a set of coupled autonomous ordinary differential equations in the form:
\begin{equation}\label{problem}
\left\{ \begin{tabular}{l}
    $\displaystyle\frac{d}{dt}{\bf x}(t) = {\bf f}({\bf x})$ \\
    ${\bf x}(t_0) = {\bf x_0}$ 
\end{tabular}  \right.
\end{equation}
where ${\mbf x}\in\mathbb{R}^m$ is the state  which depends on the time evolution $t$, ${\bf f}:\mathbb{R}^m\rightarrow \mathbb{R}^m$ is the nonlinear dynamics model, $m$ is the number of dimensions in which the problem is defined, and ${\bf x}_0$ is the initial condition of the system at time $t_0$. The KO $(\mathcal{K})$ is an infinite-dimensional linear operator that evolves all observable functions $\mathcal{G}({\bf x})$ of the state and it allows to define any problem of classical mechanics in operator form.

Let $\mathcal{F}$ be a vector space of observable functions, where $\mathcal{G}({\bf x})\in \mathcal{F}$. Since the KO is an infinite-dimensional linear operator, this space of functions $\mathcal{F}$, which the observables are defined on, is also infinite-dimensional. Therefore, if $g\subseteq\mathcal{G}({\bf x})$ is a given observable in this space, the evolution of $g$ in the dynamical system is represented by:
\begin{equation}
\mathcal{K}\left(g({\bf x})\right) = \frac{d}{dt}g({\bf x}) = \left( \nabla_{{\bf x}} g({\bf x})\right)\frac{d}{dt}{\bf x}(t) = \left( \nabla_{{\bf x}} g({\bf x})\right){\bf f}({\bf x}),
\end{equation}
where $\nabla_{{\bf x}} g = (\partial g/\partial x_1,\partial g/\partial x_2,\dots,\partial g/\partial x_m)$. That way, the evolution of any observable subjected to the dynamical system is provided by the Koopman Operator:
\begin{equation}
\mathcal{K}\left(\cdot\right) = \left( \nabla_{{\bf x}} \cdot\right){\bf f}({\bf x}).
\end{equation}

Note that the evolution of the observables is provided by the application of the chain rule for the time derivative of $g({\bf x})$. Consequently, the defined operator is linear, in that:
\begin{equation}
    \mathcal{K}\left(\beta_1g_1({\bf x})+\beta_2g_2({\bf x})\right)=\beta_1\mathcal{K}\left(g_1({\bf x})\right)+\beta_2\mathcal{K}\left(g_2({\bf x})\right),
\end{equation}
for any pair of observables $g_1\subseteq\mathcal{G}({\bf x})$ and $g_2\subseteq\mathcal{G}({\bf x})$ and any arbitrary constants $\beta_1$ and $\beta_2$. This property was outlined in Koopman's paper. The linearity of the Koopman Operator is very appealing, but this benefit is contrasted with the fact that it is infinite-dimensional. However, this issue can be overcome by capturing the evolution of the system on a finite subspace spanned by a finite set of basis functions instead of all measurement functions in a Hilbert space. In effect, this is a truncation of the Koopman operator to a finite subspace $\mathcal{F}_D$ of dimension $n$, where $\mathcal{F}_D \in \mathcal{F}$. This subspace $\mathcal{F}_D$ can be spanned by any set of eigenfunctions $\phi_i\in\mathcal{F}_D$, with $i\in\{1,2,\dots,n\}$, defined as:
\begin{equation}\label{koopman}
\mathcal{K}\left(\phi_i({\bf x})\right) =\frac{d}{dt}\phi_i({\bf x})=\lambda_i \phi_i({\bf x}),
\end{equation}
where $\lambda_i$ are the eigenvalues associated with the eigenfunctions $\phi_i$, and $n$ is the number of eigenfunctions chosen to represent the space. Therefore, the Koopman eigenfunctions can be used to form a transformation of variables that linearizes the system. Particularly, let ${\bf \Phi}({\bf x}) = \left(\phi_1({\bf x}), \dots, \phi_n({\bf x}) \right)^T$ be the set of eigenfunctions of the KO in $\mathcal{F}_D$. Then, using the relation in Eq.~\eqref{koopman}, it is possible to write the evolution of ${\bf \Phi}$ as: 
\begin{equation}\label{time_eigenf}
\mathcal{K}\left(\bf \Phi\right) = \frac{d}{dt}{\bf \Phi}=\Lambda {\bf \Phi},
\end{equation}
where $\Lambda=\text{diag}([\lambda_1, \dots,\lambda_n])$ is the diagonal matrix containing the eigenvalues of the system in $\mathcal{F}_D$. This transformation is called the Koopman Canonical Transform. The solution of Eq.~\eqref{time_eigenf} is:
\begin{equation}\label{eigen_time}
    {\bf \Phi}(t) = \exp(\Lambda t){\bf \Phi}(t_0),
\end{equation}
where ${\Phi}(t_0)$ is the value of the eigenfunctions at the initial time $t_0$. This result will be used later to solve the complete system once the eigenfunctions of the operator are obtained.

In general, we are interested in the identity observable, that is, ${\bf g}({\bf x})={\bf x}$. Therefore, it is required to be able to represent these observables in terms of the KO eigenfunctions. This is achieved using the Koopman modes, i. e., the projection of the full-state observable onto the KO eigenfunctions. If this projection can be found, the evolution of the state is represented by means of the evolution of the KO eigenfunctions, and thus, an approximate solution to the system can be provided. However, it is important to note that the challenge resides in computing the eigenfunctions and eigenvalues of the system.

For a more detailed explanation of the theory behind the Koopman eigenfunction decomposition of the dynamics, refer to Refs.~\cite{librationkoopman,zonalkoopman,analysiskoopman,servadio2021koopman}.

\section{Computing the Koopman Matrix via Galerkin Method}
This section discusses the use of the Galerkin method for computing the Eigenfunctions of the KO. First, the KO is used to define a Partial Differential Equation (PDE) for the time evolution of a scalar function $u({\bf x},t)$ (note that in general $u({\bf x},t)$ is a function of time, $t$, and ${\bf x}$). The Galerkin method is then used to convert the time evolution PDE to a matrix form using a series expansion for $u({\bf x},t)$ over a predefined basis set~\cite{librationkoopman}. This matrix form can be used to solve for the eigenfunctions and eigenvalues of the KO. 
The Koopman Operator defines a first-order PDE for the time evolution of a scale function $u({\bf x},t)$
\begin{equation}\label{koopman_PDE1}
\frac{d u({\bf x},t)}{dt}
={ f}_1\left({\bf x}\right)\frac{\partial }{\partial x_1} u({\bf x},t) +\cdots+  { f}_m\left({\bf x}\right)\frac{\partial }{\partial x_m}u({\bf x},t),
\end{equation}
Additionally, the eigenfunctions of the Koopman Operator give rise to a set of linear first-order PDEs for the eigenfunctions in the form: 
\begin{equation}\label{koopman2}
\mathcal{K}(\phi_i) = \left(\nabla_x \phi_i({\bf x})\right) {\bf f}({\bf x}) =\lambda_i \phi_i({\bf x}),
\end{equation}
or in a more expanded notation:
\begin{equation}\label{koopman_PDE}
\frac{d \phi_i({\bf x})}{dt}={ f}_1\left({\bf x}\right)\frac{\partial }{\partial x_1} \phi_i({\bf x}) +\cdots+  { f}_d\left({\bf x}\right)\frac{\partial }{\partial x_d}\phi_i({\bf x}) =\lambda_i \phi_i({\bf x}),
\end{equation}
where ${\bf f}({\bf x}) = (f_1({\bf x}), f_2({\bf x}), \dots, f_m({\bf x}))^T$. This equation is a linear first-order PDE and in general has no closed-form solution. However, it is possible to approximate the solution using the Galerkin method.

This work makes use of Legendre polynomials, due to the advantages they provide in the computation of the Koopman matrix. Legendre polynomials are a set of orthogonal polynomials defined in a Hilbert space that generate a complete basis. The idea of this methodology is to represent any function of the space by using this set of basis functions. This is done by the use of inner products and the correct normalization of the Legendre polynomials. Let $f$ and $g$ be two arbitrary functions from the Hilbert space considered. Then, the inner product between these two functions is defined as:
\begin{equation}
\langle  f, g \rangle =\int_{\Omega} f({\bf x})g({\bf x}) w({\bf x})d{\bf x} \label{innerproduct},
\end{equation}
where $w({\bf x})$ is a positive weighting function defined on the space domain $\Omega$. For the case of Legendre polynomials, the weighting function is a constant $w({\bf x}) = 1$, and the domain for each variable ranges between $[-1,1]$. In addition, the normalized Legendre polynomials are defined such that:
\begin{equation} \label{norm}
\langle \mathcal L_i, \mathcal L_j \rangle =\int_{\Omega} \mathcal L_i({\bf x}) \mathcal L_j({\bf x}) w({\bf x})d{\bf x} = \delta_{ij},
\end{equation}
where $\mathcal L_i$ and $\mathcal L_j$ with $\{i,j\}\in\{1,\dots,n\}$ are two given normalized Legendre polynomials from the set of basis functions selected, and $\delta_{ij}$ is Kronecker's delta.

\subsection{The Theory Behind the Koopman Matrix}
Any scalar function $u({\mbf x},t)$ and the KO eigenfunctions, $\phi_i({\mbf x})$, can be exactly represented as an infinite series expansion in terms of the set of basis functions. The accuracy of the Koopman approximation depends on the expansion order of the Koopman operator, which defines the number of orthogonal multivariate polynomials.
\begin{subequations}\label{series}
\begin{align}
u({\bf x},t)=\sum_{j=1}^{\infty}c_{j}(t)\mathcal L_j({\bf x})\approx \sum_{j=1}^{n}c_{j}\mathcal L_j({\mbf x})={\mbf c}^T(t)\mathcal L({\mbf x}) = \mathcal L^T({\mbf x}){\mbf c}(t)\\
\phi_i({\bf x},t)=\sum_{j=1}^{\infty}p_{ij}(t)\mathcal L_{j}({\bf x})\approx \sum_{j=1}^{n}p_{ij}\mathcal L_{\ell}({\bf x})={{\bf p}_i}^T(t){\mathcal L}({\bf x}) = {\mathcal L}^T({\bf x}){\bf p}_i(t),
\end{align}
\end{subequations}
where $c_{j}(t)$ describes the time evolution of the function $u({\bf x},t)$ over the basis $\mathcal L_{j}({\mbf x})$ and  $p_{ij}(t)$ are the coefficients associated with the eigenfunction $\phi_i$ and the basis $\mathcal L_{j}$. Moreover, ${\bf c}(t)$, ${\bf p}_{i}(t)$, and ${\mathcal L}$ are three column vectors containing the set of coefficients $c_{j}(t)$, $p_{ij}(t)$ and the whole set of basis functions, respectively. Note that although the series is infinite, a truncation was performed using $n$ different basis functions, and thus, this represents an approximation of the eigenfunctions.

The time derivative of an eigenfunction can be approximated using the Galerkin method. The general concept of Galerkin methods is to project the operator into the subspace $\mathcal{F}_D$ using a weighted residual technique such that the residual of Eq.~\eqref{koopman_PDE} is orthogonal to the span of $\mathcal{F}_D$. The eigenfunctions are computed as a linear combination of the basis functions, working directly on the coefficients $\mbf c(t)$. Starting from  $u({\bf x},t) = {\bf c}^T(t){\mathcal L}({\bf x})$, the residual error of the Koopman PDE in Eq.~\eqref{koopman_PDE1} is defined as:
\begin{equation}\label{residual0}
    e({\bf x},t) = \frac{d u ({\bf x},t)}{dt}-\mathcal{K}\left(u({\bf x},t)\right)
\end{equation}
where $\frac{d u ({\bf x},t)}{dt}$ is approximated by $\dot{u}({\bf x},t) \approx \dot{\bf c}^T(t){\mathcal L}$, with $\dot{\bf c}(t)$ a vector of size $n$. The goal is to express the derivative in the same set of basis functions as the eigenfunctions. Moreover, the solution sought should be orthogonal to $\mathcal{F}_D$, being the orthogonality condition defined as:
\begin{equation}\label{residual}
\langle \mathcal L_j({\bf x}),  e({\bf x},t)  \rangle =0, \quad \forall \ j \in \{1,2,\dots,n\}.
\end{equation}
The Koopman operator applied to $u({\bf x},t)$ is given by:
\begin{equation}\label{residual2}
\mathcal{K}\left(u({\bf x},t)\right) = \left(\nabla_{\bf x} {\mathcal L}^T{\bf c}(t)\right)^T {\bf f} =  \left({\bf f}^T\nabla_{\bf x}{\mathcal L}^T\right) {\bf c}(t)= {\bf f}^T {\mathcal L_{\mbf x}^T}{\bf c}(t),
\end{equation}
where the term ${\mathcal L_ x^T}=\nabla_{\bf x}{\mathcal L}^T$ is defined as: 
\begin{equation}\label{derivatives}
{\mathcal L_{\mbf x}^T}=\left(\begin{array}{ccc} 
    \displaystyle\frac{\partial }{\partial x_1}\mathcal L_1({\bf x}) & \cdots & \displaystyle\frac{\partial }{\partial x_m}\mathcal L_1({\bf x}) \\ 
    \vdots & \ddots & \vdots\\
    \displaystyle\frac{\partial }{\partial x_1}\mathcal L_n({\bf x}) & \cdots& \displaystyle\frac{\partial }{\partial x_m}\mathcal L_n({\bf x})\end{array}\right).
\end{equation}
Each entry in the equation can be calculated in closed-form since both ${\mathcal L}$ and ${\mathcal L_{\mbf x}^T}$ are only comprised by polynomials. Therefore, it is possible to apply Eq.~\eqref{residual2} into Eq.~\eqref{residual0} to obtain, for each basis function,
\begin{align}
    \langle \mathcal  L_j({\bf x}), e({\bf x},t)  \rangle &= \langle \mathcal L_j({\bf x}), {\mathcal L}^T({\bf x})\dot{\bf c}(t) \rangle - \langle \mathcal  L_j({\bf x}),{\bf f}^T {\mathcal L_{\mbf x}^T}{\bf c}(t) \rangle  \nonumber \\
    &=  \langle \mathcal  L_j({\bf x}), {\mathcal L}^T({\bf x}) \rangle \dot{\bf c} - \langle \mathcal  L_j({\bf x}),{\bf f}^T {\mathcal L_{\mbf x}^T} \rangle {\bf c}(t)  = 0.
\end{align}
This is a set of $n$ equations, one for each $\mathcal L_j({\bf x})$, that can be expressed in matrix form:
\begin{equation} \label{time_matrix_form}
    \frac{d {\bf c}(t)}{dt} =\mbf G^{-1} \mbf K{\bf c}(t),
\end{equation}
where:
\begin{equation} 
   \mbf  G_{ij}= \langle \mathcal  L_i({\bf x}),  \mathcal  L_j({\bf x})\rangle,
\end{equation}
However, due to the orthogonality of the basis functions we have that $\mbf G = \mbf I_{n\times n}$. Therefore, $\mbf K$ is the matrix representation of the KO over the basis $\mathcal L_i({\bf x})$, a $n\times n$ sized matrix whose components are:
\begin{equation} \label{kmatrix}
    \mbf K_{ij} = \langle \mathcal L_i({\bf x}), {\bf f}^T \nabla_{\bf x}\mathcal L_j({\bf x})\rangle.
\end{equation}
Equation \eqref{time_matrix_form}  determines the time evolution of $u({\bf x},t)$ by solving for ${\bf c}(t)$. It is a set of coupled ordinary differential equations that can be solved using the matrix exponential, ${\bf c}(t)=\exp \left(\mbf Kt\right){\bf c}(t_0)$ where ${\bf c}(t_0)$ is the initial coefficient vector.

\subsection{Evolution of Basis Functions and the Koopman Matrix}
The KO is based on the selection of a set orthogonal basis functions that describe the space of solutions. Legendre polynomials have been selected as a base for the presented application. Therefore, following the KO approach, it is desired to describe the rate of change of the basis functions, in time, as a linear combination of the functions themselves:
\begin{equation}
    \dfrac{d\mathcal L}{d t} = \tilde{\mbf K} \mathcal L \label{eqlin}
\end{equation}
where $\mathcal L$ are the Legendre polynomials and $\tilde{\mbf K} $ is the transpose of the Koopman matrix previously defined. Equation (\ref{eqlin}) is linear but high dimensional. The original ODE is expressed as a combination of linear differential equations, which solution would be exact for an infinite dimensional solution space. However, due to practicality, the Koopman solution of the ODE is an approximation that works on a subset of the infinite Hilbert space that correctly describes the solution. As such, the higher the order $c$ of the KO, the more accurate the Koopman solution, since the observables are described using more eigenfunctions that better deal with high nonlinearities. Each $i$th row of Eq.~(\ref{eqlin}) describes how each total derivative of the Legendre polynomial behaves as a linear combination of the Legendre polynomials themselves. Therefore, each entry of the Koopman matrix is evaluated through the inner product, which calculates the projection of the derivative into the basis functions. Thus, the $i$th row can be written as
\begin{align}
\dfrac{d \mathcal L_i}{d t} &= \langle \dfrac{d \mathcal L_i}{d t} , \mathcal L_0\rangle \mathcal L_0 +   \langle \dfrac{d \mathcal L_i}{d t} , \mathcal L_1\rangle  \mathcal L_1 +   \langle \dfrac{d \mathcal L_i}{d t} , \mathcal L_2\rangle  \mathcal L_2 +   \dots \nonumber \\
&= \sum_{j=0}^{n} \langle \dfrac{d \mathcal L_i}{d t} , \mathcal L_j\rangle \mathcal L_j \label{element}
\end{align}
From this expression we can write $\mbf{ \tilde K}_{ij}$ as 
\begin{equation}
\mbf{ \tilde K}_{ij}=\langle \dfrac{d \mathcal L_i}{d t} , \mathcal L_j\rangle
\end{equation}
This formulation completely describes the evolution of the system in time, for each basis function. The inner product has already been defined, and it evaluates each component of the Koopman matrix via the Galerkin method. The remaining term that requires handling is the total derivative of the basis functions. However, there is a trivial representation to best compute this term. Using chain rule, we can decompose each total derivative in terms of its partial derivatives with respect to the states. Considering the general case of a system with $m$ states, the $i$th total derivative is expanded as 
\begin{align}
\dfrac{d \mathcal L_i}{d t} &= \dfrac{\partial \mathcal L_i}{\partial x_1} \dfrac{d x_1}{d t} + \dfrac{\partial \mathcal L_i}{\partial x_2} \dfrac{d x_2}{d t} + \dfrac{\partial \mathcal L_i}{\partial x_3} \dfrac{d x_3}{d t} +  \dots \nonumber \\
&= \sum_{j=1}^{m} \dfrac{\partial \mathcal L_i}{\partial x_j} \dfrac{d x_j}{d t} \label{totd}
\end{align}
where $\mbf x$ represents the state of the system. By looking at the equation, it can be noted that each term $\dfrac{\partial \mathcal L_i}{\partial x_j}$ is always known and it does not depend directly on the problem. On the contrary, $ \dfrac{d x_j}{d t} $ is the equation of motion, which changes depending on the problem and that describes the selected application. 

Let us evolve $\mathcal L$ using Eq.~\eqref{eqlin},
\begin{equation}
    \mathcal L(t) = \exp{(\tilde{\mbf K} t)} \mathcal L (t_0)
\end{equation}
where it is assumed that $\mathcal L$ is initially orthogonal and normalized. That is, the integral
\begin{equation}
    \int_{-1}^{1} \mathcal L (t_0)\mathcal L^T (t_0) d x = \mbf I_{n\times n}
\end{equation}
where $\mbf I_{n\times n}$ is the identity matrix, respects the constraints by Eq.~\eqref{norm}. This condition has to be always valid, thus, for any time $t$, 
\begin{align}
    \int_{-1}^{1}  \mathcal L (t)\mathcal L^T (t) d x 
    &= \int_{-1}^{1} \exp{(\tilde{\mbf K} t)} \mathcal L (t_0) (\exp{(\tilde{\mbf K} t)} \mathcal L (t_0))^T d x  \\
    &= \exp{(\tilde{\mbf K} t)} \Bigg( \int_{-1}^{1} \mathcal L (t_0)  \mathcal L (t_0)^T d x \Bigg) \exp{(\tilde{\mbf K} t)}^T \\
    &= \exp{(\tilde{\mbf K} t)}  \exp{(\tilde{\mbf K} t)}^T
\end{align}
that is identity only if $\tilde{\mbf K} = - \tilde{\mbf K} ^T$, or skew-Hermitian. For this assumption, matrix $\tilde{\mbf K} $ reduces to $\mbf K $ and the KO is able to furnish the linear evolution of a function. Otherwise, further analysis that involves the Perron-Frobenius operator, which is the adjoint operator of the KO, is required \cite{Klus}.

\section{Practical Example with MATLAB - The Duffing Oscillator}
The Koopman Operator is tested in a numerical application: the Duffing Oscillator. In this chapter, we will connect the theory previously explained to its implementation in MATLAB, using the Duffing oscillator as a test case to learn how to evaluate the solution of ODEs via the KO. As such, this chapter follows and explains the MATLAB code that comes with this tutorial.

\subsection{Initial Conditions and Parameters }
The Duffing oscillator can be seen as describing the oscillations of a mass $M$ attached to a nonlinear spring,  with stiffness constant $k$, and a damper. The one-dimensional system can therefore be described by the following equations of motion, where $q$ represents the position, while $p$ is the velocity;
\begin{subequations} \label{ODE}
\begin{align} \
    \dot q &= \dfrac{p}{M}\\
    \dot p &= -k q - k a^2\epsilon q^3
\end{align}
\end{subequations}
where $a$ is a unit transformation constant and $\epsilon$ is a small parameter. Therefore, after selecting the order of the basis functions $c$ for the Koopman operator, the ODE parameters are defined.
\lstset{firstnumber=31}
\begin{lstlisting}
%% Order of the basis functions
c = 3;

%% Equation parameters
sm = 1;     % Mass
sk = 1;     % Spring constant
sa = 1;     % Unit transformation constant
se = 0.001; % Small parameter

%% Number of points in the figure
nt = 100;
\end{lstlisting}
The ODE initial condition is defined in terms of the initial state of the system $(q_0,p_0)$ and the integration final time is selected.
\lstset{firstnumber=45}
\begin{lstlisting}
q0 = 1;                 % Initial q
p0 = 0.0;               % Initial p
tf = 10;                % Final time
tk = linspace(0,tf,nt);
\end{lstlisting}

\subsection{The Differential Equation}
There are multiple different methodologies to represent a polynomial in a computer environment. In the proposed case, an approach based on a matrix representation has been selected, where different matrices express either the coefficient of each single monomial or its exponent, given a state component. 
\lstset{firstnumber=60}
\begin{lstlisting}
%% Differential equation
% Elements in the array f(a,b,c):
% a: dimension where the derivative is performed
% b: term in the derivative
% c = 1: coefficient of the term
% c \= 1: variable 'c' and its exponent in the polynomial term

% fx = [x2/sm -sk*x1-sk*se*sa^2*x1^3];
fx = zeros(2,2,3);

% Linear system
fx(1,1,1) = 1/sm; %x2/sm
fx(1,1,3) = 1;
fx(2,1,1) = -sk; %-x1*sk
fx(2,1,2) = 1;

% Perturbing term
%-sk*se*sa^2*x1^3
fx(2,2,1) = -sk*se*sa^2; 
fx(2,2,2) = 3;
\end{lstlisting}
Variable ``fx" is 3-dimensional, and it can be thought as an array (with 3 components in this particular case) of $2\times 2$ matrices. Let us analyze the entries at line 68. The first entry, ``a" in the commented code, is the number of equation of motions which matches the number of states. The second entry, ``b", is the maximum number of terms in the ODEs. Lastly, ``c" is equal to  $m+1$, where $m$ is the number of states. The first matrix, for $\text{c}=1$, stores the coefficients of each monomial of the ODE, while the remaining matrices, for $\text{c}=2,\dots,m+1$. represents the relative exponent of the state component.  Therefore, the MATLAB representation of the system of ODE from Eq.~(\ref{ODE}) is
\begin{equation}
    \text{fx(:,:,1)}=
    \begin{bmatrix}
        \sfrac{1}{M} & 0 \\
        -k & -k a^2\epsilon
    \end{bmatrix}; 
    \ \ \ \ \ \ 
    \text{fx(:,:,2)}=
    \begin{bmatrix}
        0 & 0 \\
        1 & 3
    \end{bmatrix};
    \ \ \ \ \ \ 
    \text{fx(:,:,3)}=
    \begin{bmatrix}
        1 & 0 \\
        0 & 0
    \end{bmatrix}
\end{equation}
where fx(:,:,1) are the coefficients, fx(:,:,2) are the exponents of $q$, which is the first state, and fx(:,:,3) are the exponents of $p$, which is the second and last state.

\subsection{The Basis Functions: Legendre Polynomials}
The Legendre polynomials have been chosen as the orthogonal basis functions to project the dynamics of the system. For a single variable, the polynomials $\mathcal P (\cdot)$ are:
\begin{subequations} \label{LP}
\begin{align}
    \mathcal P_0(x) &= 1 \\
    \mathcal P_1(x) &= x \\
    \mathcal P_2(x) &= \dfrac{1}{2}(3 x^2 -1) \\
    \mathcal P_3(x) &= \dfrac{1}{2}(5 x^3 -3 x) \\
    &\vdots \ \ \ \ \ \ \vdots \nonumber
\end{align}
\end{subequations}
where the subscript $i$ in $\mathcal P_i(\cdot)$ is the order of the Legendre polynomial. However, the system of polynomials expressed by Eqs.~(\ref{LP}) are mono-dimensional. The creation of basis functions in a multi-dimensional space is a problem that follows combinatorics, i.e. the theory of combinations, where we evaluate any possible combination of the polynomials among all the state components. Indeed, we can evaluate the number of basis functions for the Duffing oscillator given the KO order $c$:
\lstset{firstnumber=81}
\begin{lstlisting}
%% Number of basis functions
ns = (c + 1)*(c + 2)/2;
\end{lstlisting}
that is 
\begin{equation} \label{ns}
    n = \dfrac{(c+1)(c+2)}{2}
\end{equation}
Equation (\ref{ns}) expresses the total number of different combinations among polynomials given the maximum order. All the combinations are evaluated and stored in a matrix of dimensions $n\times m$, where row entry $i$ indicates the specific basis function $\mathcal L_i$ and each column $j$ refers to the polynomial order of the state variable $j$.
\lstset{firstnumber=84}
\begin{lstlisting}
%% Legendre basis indexes
ind = zeros(ns,2);
s = 1;
for ord = 1:c
    for i2 = 0:c
        for i1 = 0:c
            if (i1+i2) == ord
                s = s + 1;
                ind(s,:) = [i1 i2];
            end
        end
    end
end
\end{lstlisting}
For the particular case of $c=3$, the matrix that expresses the order of each single-variable polynomial into the multi-variables basis functions is
\begin{equation}
\iota = 
    \begin{bmatrix}
     0 & 1 & 0 & 2 & 1 & 0 & 3 & 2 & 1 & 0 \\
     0 & 0 & 1 & 0 & 1 & 2 & 0 & 1 & 2 & 3
    \end{bmatrix}^T
\end{equation}
The evaluation of $\iota$ can be analyzed and decomposed as 
\begin{center}
\begin{tabular}{l |c c c c c c}
 order 0 & 00 &  & & & $\rightarrow$ & 1 \\ 
 order 1 &10 & 01 & & & $\rightarrow$ & 2 \\  
 order 2 &20 & 11 & 02 & & $\rightarrow$ & 3 \\   
 order 3 &30 & 21 & 12 & 03 & $\rightarrow$ & 4 \\   
 \hline
 & &  &  &  & & $n = 10$ \\ 
\end{tabular}
\end{center}
which is the indexes analysis for the case with $m=2$ states and maximum KO order $c=3$. This provides a total of $n=10$ combinations according to Eq.~(\ref{ns}). Merging together the polynomials from Eqs.~(\ref{LP}) and the indexes from $\iota$, we can create the full array of basis functions $\mathcal L$ that expresses the Legendre polynomials for the multivariate case. Starting from order 0:
\begin{align}
    \mathcal L_0 = \mathcal P_0(q) \mathcal P_0(p) = 1
\end{align}
then order 1
\begin{align}
    \mathcal L_1 &= \mathcal P_1(q) \mathcal P_0(p) = q \\
    \mathcal L_2 &= \mathcal P_0(q) \mathcal P_1(p) = p 
\end{align}
then order 2
\begin{align}
    \mathcal L_3 &= \mathcal P_2(q) \mathcal P_0(p) =  \dfrac{1}{2}(3 q^2 -1)\\
    \mathcal L_4 &= \mathcal P_1(q) \mathcal P_1(p) = q p \\
    \mathcal L_5 &= \mathcal P_0(q) \mathcal P_2(p) =  \dfrac{1}{2}(3 p^2 -1)
\end{align}
and lastly order 3
\begin{align}
    \mathcal L_6 &= \mathcal P_3(q) \mathcal P_0(p) =  \dfrac{1}{2}(5 q^3 -3 q)\\
    \mathcal L_7 &= \mathcal P_2(q) \mathcal P_1(p) =  \dfrac{p}{2}(3 q^2 -1) \\
    \mathcal L_8 &= \mathcal P_1(q) \mathcal P_2(p) =  \dfrac{q}{2}(3 p^2 -1) \label{L8}\\
    \mathcal L_9 &= \mathcal P_0(q) \mathcal P_3(p) =  \dfrac{1}{2}(5 p^3 -3 p)
\end{align}
The set of multivariate basis functions $\mathcal L$ has been computed. However, in the current form, the multivariate Legendre polynomials are not normalized, and their inner product will produce
\begin{equation} \label{normalize}
\langle \mathcal P_i, \mathcal P_j \rangle = \dfrac{2}{2j+1}\delta_{i j}
\end{equation}
instead than just the Kronecker delta. Therefore, the polynomial normalization is necessary to respect the inner product constraint. 

The monovariable Legendre Polynomials are created in MATLAB using a recursive formulation for each coefficient of the polynomials. Given that $\mathcal P_0(x) = 1$ and $\mathcal P_1(x) = x$, each Legendre polynomial of higher order can be evaluated as 
\begin{equation}
    \mathcal P_{c+1} = \dfrac{(2 c + 1) x \mathcal P_c(x) - c \mathcal P_{c-1}(x) }{c + 1}
\end{equation}
which, in MATLAB, translates to
\lstset{firstnumber=98}
\begin{lstlisting}
%% Definition of the Legendre polynomials
LPC = zeros(c+1,c+1);
LPC(1,1) = 1;
LPC(2,2) = 1;
for i = 3:c+1
    for j = 1:i-1
        LPC(i,j+1) = LPC(i,j+1) + (2*(i-2) + 1)/(i-1)*LPC(i-1,j);
        LPC(i,j)   = LPC(i,j)   - (i-2)/(i-1)*LPC(i-2,j);
    end
end
\end{lstlisting}
Variable ``LPC" is a matrix of dimensions $(c+1)\times(c+1)$ where each row is the single-variable Legendre polynomial $\mathcal P_{i}$. More in detail, each component of the matrix is the coefficient of a monomial where the column indicates the exponential power of the variable. In the particular case of order $c=3$, matrix LPC is 
\begin{center}
\begin{tabular}{l |c c c c|}
& $x^0$ & $x^1$ & $x^2$ & $x^3$ \\
\hline
 $\mathcal P_{0}$ & 1 & 0 & 0 & 0 \\ 
 $\mathcal P_{1}$ & 0 & 1 & 0 & 0  \\  
 $\mathcal P_{2}$ & $\sfrac{-1}{2}$ & 0 & $\sfrac{3}{2}$ & 0 \\   
 $\mathcal P_{3}$ & 0 & $\sfrac{-3}{2}$ & 0 & $\sfrac{5}{2}$  \\   
 \hline
\end{tabular}
\end{center}
However, these polynomials need to be normalized such that the base is orthonormal and not only orthogonal. Therefore, each coefficient of the polynomials is normalized according to the realation in Eq.~(\ref{normalize}).
\lstset{firstnumber=109}
\begin{lstlisting}
%% Legendre polynomials in multiple dimensions
% Pre-multiplication for the normalized constant
NLPC = zeros(c+1,c+1);
for i = 1:c+1
    for j = 1:c+1
        NLPC(i,j) = sqrt((2*(i-1)+1)/2)*LPC(i,j);
    end
end
\end{lstlisting}
We created a new variable, ``NLPC", which is a matrix with the same dimension of LPC, but with normalized entries. 

The multivariate set of Legendre polynomials $\mathcal L$ is evaluated by multiplying the monovariable polynomials among each other. The order of the multiplication is dictated by matrix $\iota$, that specifies which polynomial is selected, its order, and of which state. Following the analytical computation previously performed, the MATLAB version of $\mathcal L$ is computed as
\lstset{firstnumber=118}
\begin{lstlisting}
% Multiplication of one dimensional Legendre polynomials
MLP = zeros(ns,ns);
for i = 1:ns
    for j = 1:ns
        MLP(i,j) = NLPC(ind(i,1)+1,ind(j,1)+1)*...
            NLPC(ind(i,2)+1,ind(j,2)+1);
    end
end
\end{lstlisting}
where the $i$th row of matrix ``MLP" represents the multivariate basis function $\mathcal L_i$. Matrix MLP has dimensions $n\times n$: there are $n$ rows for each Legendre polynomial $\mathcal L_i$, and $n$ columns for each possible exponent combination of the state variables in the monomials. Indeed, there is a direct connection between the rows of $\iota$ and the columns of MLP. Each row of $\iota$ specifies the exponent of the state variable for the specific coefficient. Let us have a look at the MLP matrix in the case of $c=3$:
\begin{equation}
    MLP = \begin{bmatrix}
     0.5 & 0 & 0 & 0 & 0 & 0 & 0 & 0 & 0 & 0\\  
     0 & 0.866 & 0 & 0 & 0 & 0 & 0 & 0 & 0 & 0 \\ 
     0 & 0 & 0.866 & 0 & 0 & 0 & 0 & 0 & 0 & 0\\ 
     -0.559 & 0 & 0 & 1.677 & 0 & 0 & 0 & 0 & 0 & 0\\ 
     0 & 0 & 0 & 0 & 1.5 & 0 & 0 & 0 & 0 & 0\\ 
     -0.559 & 0 & 0 & 0 & 0 & 1.677 & 0 & 0 & 0 & 0\\ 
     0 & -1.984 & 0 & 0 & 0 & 0 & 3.307 & 0 & 0 & 0\\ 
     0 & 0 & -0.968 & 0 & 0 & 0 & 0 & 2.905 & 0 & 0\\ 
     0 & -0.968 & 0 & 0 & 0 & 0 & 0 & 0 & 2.905 & 0\\ 
     0 & 0 & -1.984 & 0 & 0 & 0 & 0 & 0 & 0 & 3.307\\ 
    \end{bmatrix}
\end{equation}
If, for example, we pick the 9th row of matrix MLP, and we add information from the 9th row of $\iota$, we can reconstruct the Legendre polynomial $\mathcal L_8(q,p)$ as
\begin{equation}
    \mathcal L_8(q,p) = -0.968 q + 2.905 q p^2
\end{equation}
which is the normalized version of Eq.~(\ref{L8}).
The normalized multi-variable set of basis function $\mathcal L$ expressed in Eq.~(\ref{eqlin}) has been represented as a single two-dimensional matrix, regardless the number of state variables. All the possible combinations of monomials are ordered and stored in a second matrix, $\iota$, that works as a pointer for the whole system.

\subsection{The Koopman Matrix}
The Koopman Matrix $\mbf K$ can now be computed element by element. Looking back at Eq.~(\ref{element}), each coefficient of the Koopman matrix is evaluated through the inner product:
\begin{equation}
    \mbf K_{i j} = \langle \dfrac{d \mathcal L_i}{d t} , \mathcal L_j\rangle
\end{equation}
where the total derivative of each Legendre polynomial is evaluated according to Eq.~(\ref{totd}). While the derivatives of the states with respect to time is known (the ODE) and it is expressed by variable ``fx", the partial derivatives need to be computed. However, by being the product of different polynomials, each partial derivative of $\mathcal L_i$ can be evaluated in a trivial way:
\begin{align}
    \dfrac{\partial \mathcal L_i}{\partial x_j} &= \dfrac{\partial \left(\mathcal P_{\alpha}(x_1) \mathcal P_{\beta}(x_2)\dots \mathcal P_{\delta}(x_j) \dots \mathcal P_{\eta}(x_n) \right) }{\partial x_j} \\
    &= \left(\mathcal P_{\alpha}(x_1) \mathcal P_{\beta}(x_2)\dots \mathcal P_{\eta}(x_n) \right)\dfrac{\partial P_{\delta}(x_j)}{\partial x_j}
\end{align}
that is, the derivation affects only the mono-variable Legendre polynomial function of the state component we are deriving with respect to. For example, considering the particular case of $\mathcal L_8$ in the Duffing oscillator, from Eq.~(\ref{L8}), the partial derivatives with respect to position and velocity are
\begin{subequations}\label{partials}
\begin{align}
    \dfrac{\partial \mathcal L_8}{\partial q} &= \dfrac{\mathcal P_1(q) \mathcal P_2(p)}{\partial q} = \mathcal P_2(p) \dfrac{\partial \mathcal P_1(q)}{\partial  q} \\
    \dfrac{\partial \mathcal L_8}{\partial p} &= \dfrac{\mathcal P_1(q) \mathcal P_2(p)}{\partial p} = \mathcal P_1(q) \dfrac{\partial \mathcal P_2(p)}{\partial p} 
\end{align}
\end{subequations}
This property means that in order to evaluate the partial derivative of any multi-variable basis function $\mathcal L_i$, it is sufficient to evaluate the derivative of the mono-variable Legendre polynomial $\mathcal P_j(\cdot)$, and then multiply. The derivative of the Legendre polynomials with respect to its variable is performed as
\lstset{firstnumber=127}
\begin{lstlisting}
%% Derivative of Normalized Legendre polynomials
DLPC = zeros(c+1,c+1);
for i = 2:c+1
    for j = 1:c
        DLPC(i,j) = j*NLPC(i,j+1);
    end
end
\end{lstlisting}
where, for each $i$th polynomial, the derivatives are evaluated by multiplying the coefficients in position $j+1$ by the relative variable exponent $j$, according to the derivative law
\begin{equation}
    \dfrac{d}{d x} \left(  x^\alpha\right) =  \alpha x^{\alpha-1}
\end{equation}

The Koopman matrix $\mbf K$ is now computed inside a ``for" cycle that evaluates the coefficients row by row. Therefore, for each row of the matrix, the cycle performs three tasks. First, it evaluates the partial derivatives of the Legendre polynomials in multiple dimension, Eqs.~(\ref{partials}); then it multiplies the equations of motion (ODEs) by the partial derivatives of the Legendre polynomials to obtain the total derivative, Eq.~(\ref{totd}). Lastly, the inner product between the Legendre polynomials and the total derivative is calculated through matrix integration, Eq.~(\ref{element}). Let start with the initialization of the Koopman matrix with dimensions $n\times n$ 
\lstset{firstnumber=135}
\begin{lstlisting}
%% Operator Matrix
K = zeros(ns,ns);
MDLP = zeros(2,ns);
DB = zeros(ns*2,3);
par = zeros(1,2);
for i = 1:ns 
\end{lstlisting}
where the ``for" cycle runs along the row of $\mbf K$. The loop selects one Legendre polynomial and it evaluates the matrix coefficients one inner product at a time, following Eq.~(\ref{element}). The first step is the evaluation of the partial derivatives of the Legendre polynomials
\lstset{firstnumber=142}
\begin{lstlisting}
    % Partial Derivatives of Legendre polynomials in multiple dimensions
    for j = 1:ns
        MDLP(1,j) = DLPC(ind(i,1)+1,ind(j,1)+1)*...
            NLPC(ind(i,2)+1,ind(j,2)+1);
        MDLP(2,j) = NLPC(ind(i,1)+1,ind(j,1)+1)*...
            DLPC(ind(i,2)+1,ind(j,2)+1);
    end
\end{lstlisting}
Variable ``MDLP" is a matrix with number of rows equal to the number of states; in the Duffing oscillator $m=2$. The partial derivatives are evaluated by polynomial multiplication between the mono-variable normalized Legendre polynomials and their derivatives with respect to the states. The position of each single coefficient is dictated by the matrix of pointers $\iota$. In the loop, $i$ represents the selected Legendre polynomials we are evaluating the derivatives of.

The second step is the multiplication of the partial derivatives by the equations of motion of the system in order to obtain the total derivative with respect to time for the $i$th Legendre polynomial. Multiplication between polynomials is performed by multiplying coefficients among themselves and by adding the exponents of each variable.
\lstset{firstnumber=150}
\begin{lstlisting}
    % Total Derivative = fx * grad(Legendre Polynomials)
    DB = zeros(ns*2,3);
    s = 0;
    for dim = 1:2
        for ifx = 1:2
            if (fx(dim,ifx,1) == 0)
                break;
            else
                for j = 1:ns
                    if (MDLP(dim,j) ~= 0)
                        s = s + 1;
                        DB(s,1) = fx(dim,ifx,1)*MDLP(dim,j);
                        DB(s,2) = fx(dim,ifx,2) + ind(j,1);
                        DB(s,3) = fx(dim,ifx,3) + ind(j,2);
                    end
                end
            end
        end
    end
\end{lstlisting}
Variable ``DB" is the total derivative of the $i$th Legendre polynomial. DB is a 3 column matrix: the first column is the coefficient, the second column is the exponent of $q$, and the third column is the exponent of $p$. In the general case of a system with $m$ states, matrix DB will have $m+1$ columns. Consequently, each row of the matrix describes a monomial. The loop keeps track of the number of monomials in the polynomial thanks to the counter $s$. 

The third step for the calculation of the Koopman matrix is the evaluation of the inner product. The inner product, Eq.~(\ref{innerproduct}), consists in the integration of the product between the two entries. As such, each component of $\mbf K$ is calculated by multiplying the total derivative polynomial, DB, by the multi-variable Legendre polynomials and, afterwards, integrating in the domain of definition of the multivariate Legendre polynomials. 
\lstset{firstnumber=170}
\begin{lstlisting}
    % Matix integration
    for k = 1:ns
        for j = 1:ns
            if (MLP(k,j) ~= 0)
                for ifx = 1:s
                    flag = 1;
                    for in = 1:2
                        par(in) = round(ind(j,in) + DB(ifx,in+1) + 1);
                        if (mod(par(in),2) == 0)
                            flag = 0;
                            break;
                        end
                    end
                    if (flag == 1)
                        K(i,k) = K(i,k) + 4*MLP(k,j)*DB(ifx,1)/...
                            (par(1)*par(2));
                    end
                end
            end
        end
    end
end
\end{lstlisting}
The evaluation of the inner product is performed by working  on each single monomial and adding all the contributions together. The monomial multiplication has been divided into two different parts: exponents are summed in line 177, while coefficients are multiplied in line 184. Thus, variable ``par" stores the exponents of the monomial for each state (thus ``par" is a $m$ elements long array). Line 177 adds 1 to the exponents. This addiction correspond to the first integration step, where the exponent of the variable is increased by a unit after integration. Right after the evaluation of the exponent, an ``if" condition is checked to analyze if any exponent of the integrated monomial is even. Due to the symmetric nature of the domain with respect the origin, monomials with odd exponents before integration will result in a null contribution, since the evaluation of their integral is zero. Therefore, a flag is used to avoid the calculation of all the monomials that do not contribute to the entry of the Koopman matrix: if any of the state exponents after integration is even, then that monomial is discarded. If all the exponents of the integrated monomial are odd, the second part of the integration is performed (line 184 and 185). The monomials that are not discarded represent the integration of an even function between $-1$ and $1$. This integration can be computed in a trivial way by dividing the integration support in half and multiplying the result by a factor 2.
\begin{equation}
    \int_{-1}^{1}  \alpha x^\beta d x = 2 \int_{-0}^{1}  \alpha x^\beta d x  \ \ \ \ \ \forall \beta = 0,2,4,6,\dots
\end{equation}
In the multidimensional case, where each variable is defined in the same range, the final results is multiplied by a factor $2^m$, given $m$ the number of states. This property explains the presence of number 4 in line 184, since for the bi-variate Duffing oscillator
\begin{equation}
    \int_{-1}^{1} \int_{-1}^{1} \alpha q^\beta p^\gamma d q d p = 2^2 \int_{-0}^{1} \int_{-0}^{1} \alpha q^\beta p^\gamma d q d p \ \ \ \ \ \forall \beta = 0,2,4,6,\dots \ \wedge \ \ \forall \gamma = 0,2,4,6,\dots
\end{equation}
The remaining text in line 185 completes the integration process by diving the coefficients by the exponent of the monomial, stored in variable ``par", according to the integration law  
\begin{equation}
    \int x^\alpha d x =  \dfrac{x^{\alpha+1}}{\alpha+1} \ \ \ .
\end{equation}
A single monomial of the multiplication between the ODE (``fx") and the Legendre polynomial total derivative (``DB") has been integrated. The entry of the Koopman matrix in position $(i,k)$ is obtained by adding the contributions of all the monomials: this explains the recursive formulation of line 184, where $\mbf K_{i,j}$ redefines itself after each single iteration. The Koopman matrix is therefore calculated by filling each column, $k$, and then moving to the next row, $i$.

\subsection{The Selection of the Observables}
The system is described through the Koopman matrix. However, in order to obtain the solution of the system, we need to obtain the state variables as a linear combination of the basis functions. Indeed, thanks to the basis functions, any function of the state can be approximated as a linear combination of the basis functions. Let $\mathbf g (\mbf x)$ be a set of observables that we are interested to evaluate. Any $i$th observable can be projected on the set of basis functions using the inner product:
\begin{equation}
    \mbf g = \sum_{j=0}^{n} \langle \mbf g, \mathcal L_j\rangle \mathcal L_j  
\end{equation}
These observables can be any function of the original variables $\mbf x$, including the states themselves, called the identity observables. Therefore, the set of coefficients that expresses the projection of the observables into the basis functions can be represented in matrix form as in the case of the Koopman matrix:
\begin{equation}
    \mathcal H_{i,j} = \langle \mbf g_i , \mathcal L_j\rangle
\end{equation}
where $\mathcal H_{i,j}$ indicates the projection of the $i$th observable onto the $j$th basis function. Any number of observables can be evaluated by increasing the number of rows of $\mathcal H$. In the particular case of the identity observables, where it is desired to obtain the state of the system, matrix $\mathcal H$ has dimensions $m \times n$. 

The evaluation of the observable matrix is analogous of the calculation of the Koopman matrix presented in the previous chapter. This time, the two polynomials multiplied in the inner product are the Legendre polynomials and the observable polynomials. 
\lstset{firstnumber=193}
\begin{lstlisting}
H = zeros(2,ns);
Obs = zeros(ns*2,3);
for i = 1:2
    
    % Observable Polynomials
    Obs = zeros(ns*2,3);
    Obs(1,1) = 1;
    if (i==1)
        Obs(1,2) = 1; %q
    elseif (i==2)
        Obs(1,3) = 1; %p
    end
\end{lstlisting}
The observable polynomials are represented with the same notation used for the total derivatives. In the first iteration, the position polynomial is created, which consists in the single monomial $1q^1 \rightarrow \begin{bmatrix} 1 & 1 & 0\end{bmatrix}$. For the second iteration, the velocity is selected, and thus the observable matrix representation of the polynomials in the code is $1p^1 \rightarrow \begin{bmatrix} 1 & 0 & 1\end{bmatrix}$. Each entry of $\mathcal H$ is evaluated through polynomial integration, likewise for $\mbf K$, according to the inner product.
\lstset{firstnumber=207}
\begin{lstlisting}
% Matix integration
    for k = 1:ns
        for j = 1:ns
            if (MLP(k,j) ~= 0)
                flag = 1;
                for in = 1:2
                    par(in) = round(ind(j,in) + Obs(1,in+1) + 1);
                    if (mod(par(in),2) == 0)
                        flag = 0;
                        break;
                    end
                end
                if (flag == 1)
                    H(i,k) = H(i,k) + 4*MLP(k,j)*Obs(1,1)/...
                        (par(1)*par(2));
                end
            end
        end
    end
end
\end{lstlisting}

\subsection{Eigendecomposition of the System }
In the previous parts of the algorithm, the Koopman matrix is calculated, the basis function are defined, and the nonlinear equations of motion have been approximated as a finite-dimensional linear differential equation in the form of Eq.~(\ref{eqlin}). Under the assumption of matrix $\mbf K $ being diagonalizable, the eigendecomposition of the dynamics can be performed;  
\begin{equation}
    V \mbf K = \Lambda V   \label{eig}
\end{equation}
where $V$ is the eigenvectors matrix and $\Lambda$ is the diagonal matrix of eigenvalues. For high Koopman operator order, matrix $\mbf K$ could be ill conditioned. Therefore, in MATLAB, it is convenient to add the ``nobalance" option during diagonalization. 
\lstset{firstnumber=228}
\begin{lstlisting}
%% Eigenvalue decomposition
[V,D] = eig(K,'nobalance');
iV = inv(V);
\end{lstlisting}
If instead, the resultant matrix $K$ is non-diagonalizable, a different approach is required such as the one presented in Ref.~\cite{schur} which covers both diagonalizable and non-diagonalizable systems.

\subsection{The Theory Behind Diagonalization and the Time Solution of the Observables}
Before continuing with the code and evaluating the time evolution of the state of the system, it is important to take a step back and analyze how the solution is obtained. Let us rewrite Eq.~(\ref{eqlin}) by highlighting the variable dependency of the functions:
\begin{equation}
    \dfrac{d}{d t} (\mathcal L(\mbf x(t)))= \mbf K \mathcal L(\mbf x(t)) \label{eqlin_dep}
\end{equation}
Thanks to the eigendecomposition of the Koopman matrix, Eq.~(\ref{eig}), we can write the eigenfunctions of the system from the Legendre polynomials using the eigenvectors matrix:
\begin{equation}
    \boldsymbol \phi(\mbf x(t)) = V \mathcal L(\mbf x(t)) \label{eigfun}
\end{equation}
The analysis of the eigenfunctions brings us to a simpler description of the system where each differential equation is decoupled from the others, and, therefore, easier to solve. Let us take the time derivative of the eigenfunctions: 
\begin{subequations}
\begin{align}
    \dfrac{d}{d t }\boldsymbol \phi(\mbf x(t)) &= \dfrac{d}{d t }  (V \mathcal L(\mbf x(t))) = V \dfrac{d}{d t }  ( \mathcal L(\mbf x(t))) \\
     &= V \mbf K \mathcal L(\mbf x(t)) \\
     &= \Lambda V \mathcal L(\mbf x(t)) \\
     &= \Lambda \boldsymbol \phi(\mbf x(t)) \label{diag}
\end{align}
\end{subequations}
where the substitutions come, respectively, from Eq.~(\ref{eqlin_dep}), Eq.~(\ref{eig}), and Eq.~(\ref{eigfun}). The last relation, Eq.~(\ref{diag}), shows a diagonal system of ODEs whose solution is known
\begin{equation}
    \boldsymbol \phi(\mbf x(t)) = \exp(\Lambda t) \boldsymbol \phi(\mbf x(t_0)) \label{sol}
\end{equation}
where $\boldsymbol \phi(\mbf x(t_0))$ indicates the value of the eigenfunctions at the initial time $t_0$. 

The evolution of the eigenfunctions with time is used to find the solution of any observable, and thus, of the state of the system (when the observable is the identity). Each observable function $\mbf g (\mbf x (t))$ has already been represented as a linear combination of the basis functions of the system through matrix $\mathcal H$. Therefore, after some manipulations, the time evolution of the observables can be expressed solely as a function of time:
\begin{subequations}
\begin{align}
   \mbf g (\mbf x (t)) &= \mathcal H \mathcal L(\mbf x(t)) \\
   &= \mathcal H V^{-1} \boldsymbol \phi(\mbf x(t))\\
   &= \mathcal H V^{-1} \exp(\Lambda t) \boldsymbol \phi(\mbf x(t_0))\\
   &= \mathcal H V^{-1} \exp(\Lambda t) V \mathcal L(\mbf x(t_0))  \label{final}
\end{align}
\end{subequations}
where the substitutions come, respectively, from the inversion of Eq.~(\ref{eigfun}), Eq.~(\ref{sol}), and Eq.~(\ref{eigfun}). The solution of the equations of motion has been found, having picked the state of the system as observables. 

Few considerations and comments are due before continuing with the remaining part of the code. Equation~(\ref{final}) has a single dependency on time only, meaning that the value of the observable at any give time-step can be known by mere evaluation of the function. Moreover, matrix $\mathcal H$ is the only term that directly depends on the selection of the observable: the remaining part of the formulation describes the time evolution of the eigenfunctions (and thus, basis functions). That is, once the Koopman matrix is available, a different observable can be evaluated just by calculating a new matrix $\mathcal H$. The observable matrix specifies the coefficients in front of the eigenfunctions, such that any observable is a different linear combination of the latter.

\subsection{The Eigenfunctions}
The Koopman solution is evaluated once given an initial condition of the state of the system. This initial condition is transformed into the initial values of the basis functions at time $t_0$. Therefore, the Legendre polynomials are evaluated at a given state, providing $\mathcal L(\mbf x(t_0))$.
\lstset{firstnumber=232}
\begin{lstlisting}
%% Generation of the functional space of solutions
PHI = H * V;

phi0 = zeros(ns,1);
h0 = zeros(ns,1);
for i = 1:ns
    for j = 1:ns
        h0(i) = h0(i) + MLP(i,j)*q0^ind(j,1)*p0^ind(j,2);
    end
end
\end{lstlisting}
The loop is nothing more than a mere polynomial evaluation, where we substitute the numerical outcome from the variables of the polynomials themselves. Once again, $\iota$ specifies the exponent of the variables for each monomial, while variable ``MLP" provides the coefficients. The initial condition in terms of the eigenfunctions of the system is evaluated by projecting the Legendre polynomials using the eigenvectors.
\lstset{firstnumber=243}
\begin{lstlisting}
% Initial Conditions
for i = 1:ns
    phi0(i) = iV(i,:)*h0;
end
\end{lstlisting}

\subsection{The Koopman Solution}
Finally, everything is set to provide the final Koopman solution of the state of the system. Thus, the numerical values of the state for any given time step is provided by solving Eq.~(\ref{final}) and extracting the real part:
\lstset{firstnumber=248}
\begin{lstlisting}
%% Computaion of the solution as a function of time
Sol = zeros(2,nt);
for k = 1:nt
    Sol(:,k) = real(PHI*diag(exp((tk(k)*diag(D))))*phi0);
end
q = Sol(1,:);
p = Sol(2,:);
\end{lstlisting}
where each row of the solution matrix ``Sol" correspond to the time behaviour of a single state. 

The accuracy of the solution improves by increasing the number of eigenfunctions that represent the state space, thus by augmenting the KO order $c$.

\section{Conclusions}
The state of the system has been expressed as a linear combination of the eigenfunctions. The value of position and velocity at any time step can be directly evaluated by a simple polynomial evaluation of the basis functions. Consequently, the Koopman approach offers an analytical analysis of the system, where the dynamics have been projected onto a well-defined set of orthogonal functions. The accuracy of the Koopman approximation of the solution is improved by increasing the number of eigenfunctions, and thus, by increasing the Koopman order $c$. Therefore, the initial nonlinear system is now expressed with a new, linear, highly-dimensional representation. The reader is advised to compare the KO solution with the numerical integration offered by MATLAB ODE45 function (comparison included in the code that comes with this tutorial). The results show the dependency of the solution on the Koopman order and on the intensity of the small parameter $\epsilon$.

\section{Citation of this Tutorial}
The aim if this tutorial is to encourage the use of the Koopman Operator among different disciplines and to share a new, reliable, methodology to evaluate the Koopman matrix. If this tutorial proves to be helpful to the reader, the authors kindly ask for this tutorial to be cited by any works that it inspires.

\end{document}